\documentclass{amsart}
\usepackage{graphicx} 
\usepackage[a4paper,left=2.25cm,right=2.25cm,top=3cm,bottom=3cm]{geometry}

\usepackage{
  amsmath,
  amssymb,
  amsthm,
  paralist,
  thmtools,
  tikz-cd,
  verbatim,
  relsize, 
  bbm, 
  dsfont,
  mathtools,
  tikz-3dplot, 
  tikz,
  xcolor,
  multirow
}
\usepackage{hyperref}
\usepackage{soul}

\newcommand\ZZ{{\mathbb{Z}}}

\title{Generating Hadamard matrices with transformers}
\date{April 2026}

\author{Geordie Williamson}
\address{Geordie Williamson, School of Mathematics and Statistics and the Sydney Mathematical Research Institute, The University of Sydney, 
Camperdown 2005, Australia}
\email{g.williamson@sydney.edu.au}
\thanks{GW is supported by ARC grant DP230102982 and the Max Planck Humboldt Research Award. He is grateful for the use of the National Computational Infrastructure (NCI) during this project.}
\author{Oded Yacobi}
\address{Oded Yacobi, School of Mathematics and Statistics, The University of Sydney, 
Camperdown 2005, Australia}
\email{oded.yacobi@sydney.edu.au}
\thanks{OY is supported by ARC grant DP230100654. OY thanks Dmitrii V. Pasechnik for useful discussions on databases of Hadamard matrices.}
\author{Paul Zinn-Justin}
\address{Paul Zinn-Justin, School of Mathematics and Statistics, The University of Melbourne, 
Victoria 3010, Australia}
\email{pzinn@unimelb.edu.au}
\thanks{PZJ 
acknowledges the use of the computing facilities of the University of Melbourne
(High Performance Computer system and Research Computing Portal).}

\begin{document}

\maketitle

\begin{abstract}

We present a new method for constructing Hadamard matrices that combines transformer neural networks with local search in the PatternBoost framework. Our approach is designed for extremely sparse combinatorial search problems and is particularly effective for Hadamard matrices of Goethals--Seidel type, where Fourier methods permit fast scoring and optimisation. For orders between $100$ and $200$, it produces large numbers of inequivalent Hadamard matrices, and for larger orders, it succeeds where local search from random initialisation fails. The largest example found by our method has order $252$. In addition to these new constructions, our experiments reveal that the transformer can discover and exploit useful hidden symmetry in the search space.

\end{abstract}

\section{Introduction}

\subsection{Background} A Hadamard matrix is an $n\times n$ matrix with entries $\pm1$ whose rows are mutually orthogonal. 
A  necessary condition is that, apart from the trivial orders $1$ and $2$, the order $n$ must be a multiple of $4$. Hadamard’s conjecture is that this necessary condition is also sufficient. The problem goes back to Hadamard’s 1893 work on the maximal determinant problem \cite{hadamard1893resolution}: a matrix $M$ is Hadamard if and only if its entries have absolute value $\leq 1$ and $|\det(M)|$ achieves the maximal value $n^{\frac{n}{2}}$.

Hadamard's Conjecture is considered one of the central open problems in combinatorial design theory, in part because Hadamard matrices have  practical applications in signal processing and related areas \cite{yarlagadda2012hadamard}. 
Currently, the only orders less than 1000 for which no Hadamard matrix is known to exist are 668, 716, 892. The last breakthrough was the construction of a Hadamard matrix of order 428 in 2005 \cite{H428}.

Even in ranges where Hadamard matrices are known to exist, the data can be quite sparse. For instance, Sloane's database of Hadamard matrices \cite{Sloane_Hadamard} contains a single example in each possible order up to $256$, and for the highest values of $n$, these appear to be the only known examples (up to equivalence) in the literature. See also \cite{CatiPasechnik2024} for an overview of known Hadamard constructions.

The aim of this paper is three-fold. Firstly, we generate many new examples of Hadamard matrices for orders $100 \leq n \leq 252$.
In fact, our construction is inherently probabilistic, and we can produce new Hadamard matrices every time we run the construction. 
For example, we can produce around 5\,000 $n=148$ Hadamard matrices per hour on a NVIDIA H100 GPU.

Secondly, we use transformers in our constructions, and we believe the problem of generating many or large Hadamard matrices is an excellent benchmark in the emerging field of ``AI for Math''. Indeed, this problem is easily stated and implemented, has instances of greatly varying difficulty, and has historically led to several rich veins of mathematical development. We hope that our paper may inspire other AI-approaches to this challenging problem.\footnote{We are not the only ones to think so, see \cite{Epoch_Hadamard}.} Existing attempts to use neural networks on this problem include \cite{peres2022equivariantneuralnetworksrecovery,dhanaraju2023deep}.

Finally, we hope our approach shows that transformer-guided search can be an effective tool for combinatorial design, and suggest a broader role for AI-based methods in difficult existence and classification problems in combinatorics.

\subsection{The role of transformers} Transformers are a neural network architecture designed to process sequences using an attention mechanism. 
We use transformers in a  method called PatternBoost \cite{CEWW}, which  is a variation on classical evolutionary algorithms, designed to find a ``needle in a haystack''.\footnote{For an excellent and flexible implementation of PatternBoost, see \cite{axiommath2025axplorer}.}


PatternBoost applies to search problems where one can randomly generate many samples, and score them so that samples with better scores are ``nearer'' to the desired output. For example, in our situation we can easily generate $\pm1$ matrices and score them by the function $n^{\frac{n}{2}}-|\det(M)|$\footnote{We use a more sophisticated score function in practice, see \S\ref{ssec:score}.}.
The PatternBoost implementation proceeds by randomly generating a large sample, culling the sample using local search techniques and then training a neural network on the result. This essentially produces a probability distribution, from which we can again generate a large sample and repeat. This process is depicted below:
\vspace{2mm}

\begin{center}
\begin{tikzpicture}[node distance=10ex]
\tikzstyle{object} = [rectangle, text centered, draw=black, dotted, font=\footnotesize, fill=yellow!30]
\tikzstyle{field} = [rectangle, text centered, draw=black, dotted, font=\footnotesize, fill=blue!10]
\tikzstyle{diagram} = [thick, fill=violet!20, rounded corners, text centered, draw=black, font=\footnotesize]
\tikzstyle{answer} = [fill=green!20, rounded corners, text centered, dotted, draw=black, font=\footnotesize]

\tikzstyle{implies} = [double,double equal sign distance,-implies, thick]
\tikzstyle{arrow} = [dashed, arrows = {-Stealth[length=10pt, inset=6pt]}]
\tikzstyle{projectarrow} = [thick, arrows = {-Stealth[length=10pt, inset=6pt]}]

\node (LS) [diagram] {large sample};
\node (IS) [diagram, right of = LS, xshift=20ex] {improved sample};
\node (TM) [diagram, right of = IS, xshift=20ex] {trained model};

\draw[projectarrow] (LS) to node[midway, above, font=\footnotesize]{local search} (IS);
\draw[projectarrow] (IS) to node[midway, above, font=\footnotesize]{transformer} (TM);
\draw[projectarrow] (TM) to[bend left=17] node[midway, above, xshift=-1ex, font=\footnotesize]{sample from model} (LS);
\draw[->, thick, >=Stealth] (-4.5,0) -- (-1.1,0) node[midway, above] {\footnotesize{randomly generate}};
\end{tikzpicture}
\end{center}

Note that the reproduction step uses transformers, and the hope is that this produces samples that more effectively inherit favorable traits from previous generations as compared to more classical reproduction algorithms.

\subsection{The role of symmetry}\label{ssec:sym}
Let $\mathcal{M}_n$ denote the set of $n\times n$ matrices with $\pm1$ entries, and  
let $G=C_2^n \rtimes S_n$ denote the group of signed permutations. There is an obvious action of $G\times G$ on $\mathcal{M}_n$ by left and right multiplication, which preserves the subset of Hadamard matrices. Two Hadamard matrices are \emph{equivalent} if they are in the same $G\times G$-orbit. One typically considers score functions that are invariant under $G\times G$.

This symmetry suggests restricting attention to structured representatives within each equivalence class, particularly those exhibiting additional internal symmetries. We will therefore be interested in particular forms of Hadamard matrices that possess a large symmetry.
Williamson \cite{Williamson} first suggested the use of a particular “quaternionic” form to construct Hadamard matrices, namely
\begin{equation}\label{eq:will}
M = \begin{bmatrix}
    A & B & C & D
    \\
    -B & A & -D & C
    \\
    -C & D & A & -B
    \\
    -D & -C & B & A
\end{bmatrix}
\end{equation}
where $A,B,C,D$ are symmetric circulant matrices of size $n'=n/4$.

\begin{figure}[h!]
    \includegraphics{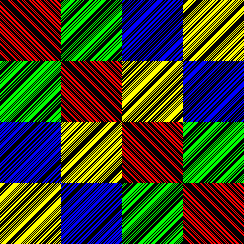}
    \caption{A GS type Hadamard matrix of order $244$. Coloured entries are $+1$, black entries are $-1$.}
    \label{fig:H244}
\end{figure}

Among the many variations of this construction, we shall focus mostly on
the Goethals--Seidel \cite{GS} form
\begin{equation}\label{eq:GS}
M = \begin{bmatrix}
    A & BF & CF & DF
    \\
    -BF & A & -FD & FC
    \\
    -CF & FD & A & -FB
    \\
    -DF & -FC & FB & A
\end{bmatrix}
\end{equation}
where $A,B,C,D$ are circulant matrices, and $F$ is the square $n'\times n'$ matrix with $1$s on the antidiagonal and $0$s elsewhere. Note that a matrix $M$ as in \eqref{eq:GS}, which we will call of ``GS type", is completely determined by its first row.
Fig.~\ref{fig:H244} provides an example of a Hadamard matrix of GS type, as generated by our program.

In both constructions, $\det M$ is easily computed in terms of the eigenvalues of $A,B,C,D$, i.e., the discrete Fourier modes of one of their rows.

\begin{figure}[h!]
\resizebox{\linewidth}{!}{\input{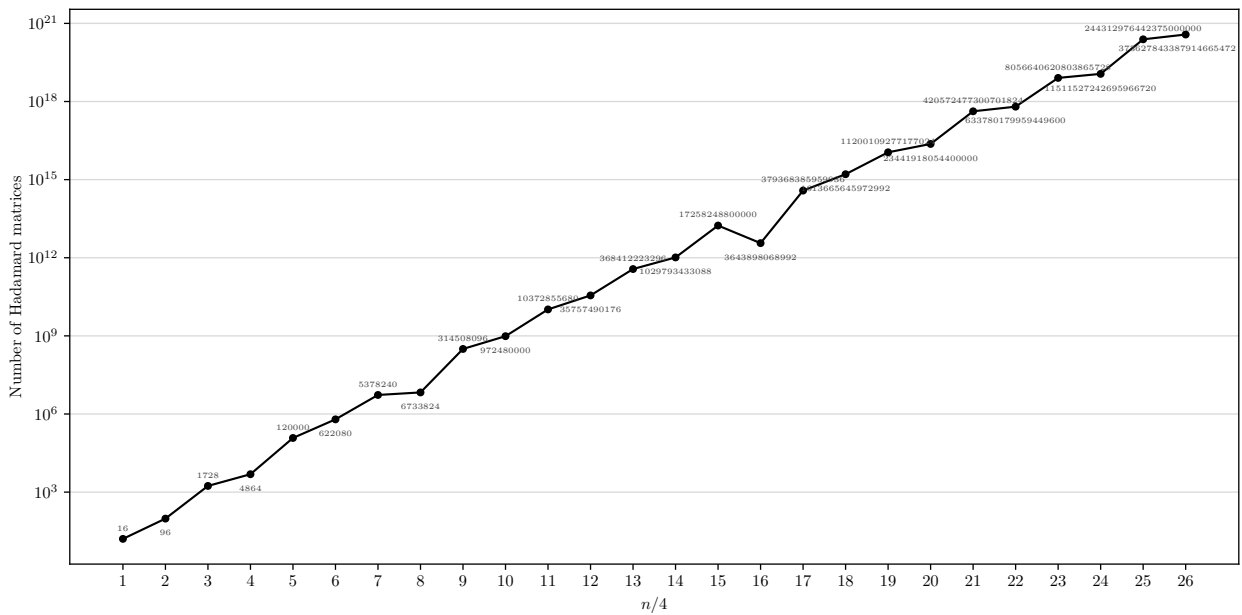}}
\caption{Number of GS type Hadamard matrices for small $n$.}
\label{fig:exact_GS}
\end{figure}

One reason to use matrices of GS type is that we expect Hadamard matrices of this type to be ``abundant'': based on small $n$ exact enumeration (see Figure~\ref{fig:exact_GS}),
we expect their number to grow like $\approx 1.56^n$, and so dividing by the total number of matrices of GS type, we obtain a ratio of $\approx 0.78^n$. For comparison, one expects that the ratio of {\em all}\/ Hadamard matrices among $\{\pm1\}$-valued matrices to decrease like $\alpha^{n^2}$.
There is a residual symmetry group acting on GS type matrices,
namely $H:=(D_{2n'}\times \ZZ/2\ZZ)^4 \rtimes \mathcal S_4 \times \text{Aut}(\ZZ/n'\ZZ)$. Although $H$ is much smaller than $G\times G$, it will play a key role in our algorithm (see in particular \S\ref{ssec:transf}).

\subsection*{Code and data}
All our code is available on GitHub at \url{https://github.com/pzinn/hadamard};
many runs, including all the ones that were used in the making of this paper,
are recorded at \url{https://wandb.ai/aiformath/hadamard}.

\section{Algorithms and implementation details}

The PatternBoost algorithm is iterative and proceeds as follows.
At generation $0$ we start from a large batch of random samples of arrays encoding matrices of the prescribed (GS) type. Each batch is then improved by a collection of local and non-local optimisation procedures, and scored. The resulting arrays are deduplicated after reducing to a canonical form under the symmetry group $H$.
We keep only the best part of this pool, both as the current set of candidates and as the training data for the transformer update. The next generations follow the same principle, except new samples are produced by the transformer rather than by random initialisation, and merged with the previous population.
Note that the transformer is not asked to solve the problem on its own: it learns from the output of the optimisation stage, while the optimisation stage in turn starts from transformer-generated proposals. This feedback loop is the basic PatternBoost mechanism used throughout the paper. 

The following subsections describe the main ingredients of this loop: the transformer model, the score function, the use of segment sums to reduce the search space, the improvement procedures, and some remarks on implementation.

\subsection{Transformer}\label{ssec:transf}
Our basic model is a decoder-only transformer trained autoregressively on arrays encoding the first rows of the circulant blocks. The array is first packed into ``stacks'' of $s$ bits, where $s$ is the \verb|stacking| parameter; these stacks are treated as tokens in a vocabulary of size $2^s$. The main architectural hyperparameters, besides the stacking parameter, are the number of layers, the embedding dimension, and the number of attention heads (in the code,
\verb|n_layer|, \verb|n_embd|, \verb|n_head|). 

Training is standard next-token prediction with cross-entropy loss, but with an important modification: the training examples are randomly transformed on the fly by the symmetry group of the search space (GS type matrices) $H=(D_{2n'}\times \ZZ/2\ZZ)^4 \rtimes \mathcal S_4 \times \text{Aut}(\ZZ/n'\ZZ)$. Concretely, the code applies random cyclic shifts, reversals, sign changes, segment permutations, and automorphisms of the underlying cyclic group. 
This serves two purposes. First, it prevents the model from locking onto arbitrary representatives of a symmetry class. Second, because of how large the symmetry group $H$ is,\footnote{For instance, if $n'$ is prime, and dividing by $\ZZ_2$ due to overcounting the reversal, the symmetry group is of order $\sim \frac{3}{4}n^5$.} this procedure makes the effective training set essentially inexhaustible, so that overfitting is never observed in practice, eliminating the need for a test set. 

The model is trained with the AdamW optimiser and cross-entropy loss. At generation~$0$ a full training run is performed; at subsequent generations, training resumes from the previous model with a reduced number of steps and learning rate, since the distribution shifts only incrementally between generations.
The main training hyperparameters are simply the \verb|training_size|,
the \verb|learning_rate|, and the (number of) \verb|training_steps|.

During the sampling phase, the transformer 
generates autoregressively a population whose size is controlled by the parameter \verb|sample_size|.
(Typically, \verb|training_size| is a small fraction of \verb|sample_size| -- 5\% or 10\% in most experiments.)
A \verb|temperature| parameter controls the trade-off between concentrating on the most likely continuations and preserving diversity. In our experiments, a moderate temperature works best: too small a value produces many near-duplicates, while too large a value destroys the structural patterns learned during training. 
A common heuristic is to increase the temperature from one generation to the next, so that the model first exploits what it has already learned and only later injects more diversity. 

\subsection{Score function}\label{ssec:score}
Let $S=-\log\det(M/n^{1/2})$;
it is nonnegative for all $\{\pm1\}$-valued (in fact,
$[-1,1]$-valued) matrices, and zero iff $M$ is Hadamard.

For matrices of GS type \eqref{eq:GS}, with the first rows of $A,B,C,D$ given by $a_{i,j}$, $i=1,\ldots,4$, $j=0,\ldots,n'-1$,
the score function can be rewritten
\begin{equation}\label{eq:score_fft}
S=-\log\det(M/n^{1/2}) = -\sum_{j=0}^{n'-1} \log\left(\sum_{i=1}^4 |\hat a_{i,j}|^2\right)
\end{equation}
in terms of the Fourier modes $$\hat a_{i,j} = \frac{1}{\sqrt{n}}\sum_{k=0}^{n'-1}a_{i,k}\omega^{jk}, \omega=\exp(2\pi i/n').$$
Its implementation is particularly fast (FFT), and its variation under local modification can also be computed quite efficiently.

Note that because of the sum rule 
$\sum_{j=0}^{n'-1}\sum_{i=1}^4 |\hat a_{i,j}|^2=1$, $S$ is minimum when {\em each term}\/ in the summation over $j$ is zero. For implementation purposes, one should make the minimum more explicit by writing
$$S=\sum_{j=0}^{n'-1} f\left(\sum_{i=1}^4 |\hat a_{i,j}|^2-1\right),$$
with $f(u)=u-\log (1+u)$ to minimise rounding errors. In principle $f$ can be replaced
with any function with a single minimum at $u=0$, although we have found that the 
above choice works best for the purpose of local search.

\subsection{Segment sums}\label{ssec:ss}
The zero Fourier modes provide us with a well-known constraint on the ``segment sums'' $k_i:=\sqrt{n}\,\hat a_{i,0}=\sum_{j=0}^{n'-1} a_{i,j}$, namely
\begin{equation}\label{eq:ss}
\sum_{i=1}^4 k_i^2 = n
\end{equation}
This equation admits only finitely many solutions for $k_i$ integers of the same parity as $n'$. The code allows to fix these numbers (parameter
\verb|segment_sums|) to reduce the search space (asymptotically,
it reduces it by a factor of $(8/(\pi e n))^2\approx .877589 / n^2$).
As we shall discuss next, imposing this constraint has significant impact on local search algorithms. The design of the transformer is not affected by the value of \verb|segment_sums|; however, if trained on a data set with a fixed choice of segment sums, the model quickly learns to sample arrays that satisfy this constraint with very high probability. Further discussion of this feature is given in \S\ref{ssec:discuss_ss}.

\subsection{Improvement}\label{ssec:improv}
The improvement stage combines several heuristics. The first one is an exhaustive one-bit descent, implemented efficiently through incremental Fourier updates. Starting from the best one-bit moves, the code also tests selected multi-bit combinations, which helps escape the obvious local minima of purely greedy descent. This is complemented by random non-local moves, one segment at a time:
fixing the Fourier transforms of the other three segments, the target Fourier transform for the remaining segment is computed (the one that would make the matrix Hadamard), and the block is replaced by the closest $\{\pm1\}$-vector to the inverse Fourier transform of a random phase perturbation of this target. This can make large jumps in the search space and is particularly effective when the array is already close to a Hadamard matrix.
The parameter {\tt num\_improve} controls how aggressively these improvement routines are applied between two successive transformer updates. In particular, if {\tt num\_improve}$>0$,
a parallel tempering step is first applied to the data, with the Metropolis algorithm; it runs several replicas at different effective temperatures and occasionally swaps them, allowing difficult configurations to cross energetic barriers that would trap a purely gradient descent search. The temperature ladder is autotuned during the run to maintain swap acceptance rates in a target range.

The case of fixed segment sums deserves special treatment: one can no longer perform the most efficient one-bit-flip searches, because such moves leave the constrained search space. As a result, imposing segment sums shrinks the search space but also makes local search less efficient, and the code must rely on more indirect admissible moves (such as cyclic rotations of entries) within the constrained family.

\subsection{Modifying the transformer}\label{ssec:mod_transf}
A natural design question is whether one should stay with a completely generic transformer, or instead inject more problem-specific information into the model. The code supports both philosophies. The vanilla version is simpler and more portable: it treats the encoded array as an ordinary token sequence and leaves all structural inference to training. The modified version is more specialised and attempts to expose additional information about partial constructions to the model. This makes the architecture less generic, but may improve sample efficiency in a search problem where the target set is extremely sparse.

The main such modification is controlled by the flag \verb|transformer_uses_score|. In this mode, the four circulant blocks are generated one after another rather than as a single long token sequence, and the model is given a coarse summary of how well the already generated blocks fit the Fourier constraints. This reduces the context window from the full order $n$ to the segment length $n'$ and lets the model condition its later choices on partial score information. Conceptually, the vanilla transformer learns only from previously generated symbols, whereas the modified transformer is also informed by a low-dimensional state variable measuring partial progress toward the Hadamard conditions. 

A key design choice is that training and inference use slightly different score contexts. During training, the model sees the actual cumulative spectrum of the training example (which typically exceeds the Hadamard target of~$1$ per frequency). During inference, the process starts from the ideal target spectrum. This deliberate mismatch encourages the model to extrapolate toward the Hadamard condition rather than simply reproduce the distribution of the training data.

See \S\ref{ssec:discuss_mod} for a comparison of the performance of the vanilla and modified transformers.

\subsection{Implementation}
The implementation is written in PyTorch and runs on a single GPU. 
All operations -- scoring, local search, parallel tempering, transformer training and sampling -- are fully batched and parallelised over the candidate arrays. The local search routines maintain the Fourier transforms of all arrays simultaneously and use batched tensor operations to evaluate and apply moves across the entire population at once. 
This is essential because the algorithm relies on generating and filtering very large populations of candidates. 
Symmetry considerations are built into the algorithm; in particular,
the same symmetry operations serve both for data augmentation during training and for canonicalisation and deduplication of solutions.

A typical run at $n=188$ with a sample size of $1\,000\,000$ and $30$ generations takes a few hours on an NVIDIA H100 GPU.

\section{Results of experiments}
The experiments that are reported below were kept as simple as possible, so without some of the bells and whistles of the code (in particular, the \verb|temperature| is kept equal to $1$ for all sampling); the hyperparameters are given typical reasonable values -- we have not tried systematic grid searches for optimal parameters.
The size of the transformer model is modest by modern standards, though, as we shall see, it is large enough to perceive some nontrivial structure in the data it is trained on and sample adequate Hadamard candidates.

\subsection{Mass production of Hadamard matrices}
\begin{figure}[h!]
\def\runname{aehw4wce}
    \includegraphics[width=10cm]{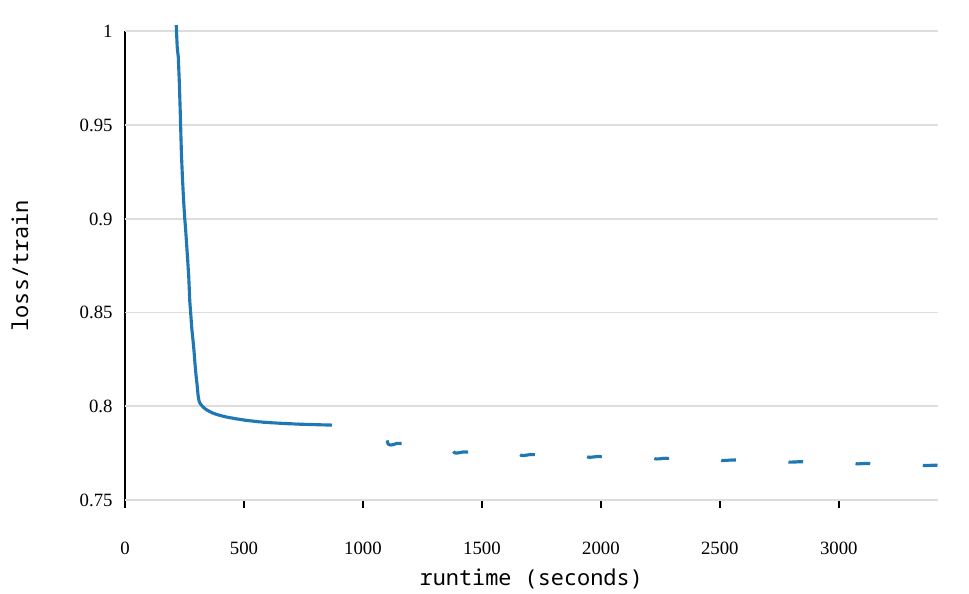}
    \\[2mm]
\includegraphics[width=7cm]{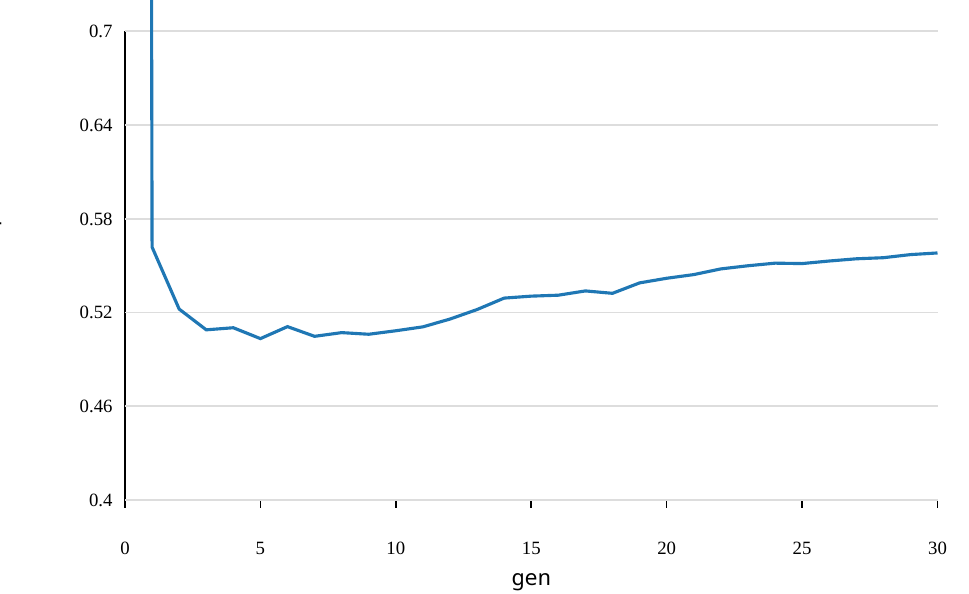}
\qquad
\includegraphics[width=7cm]{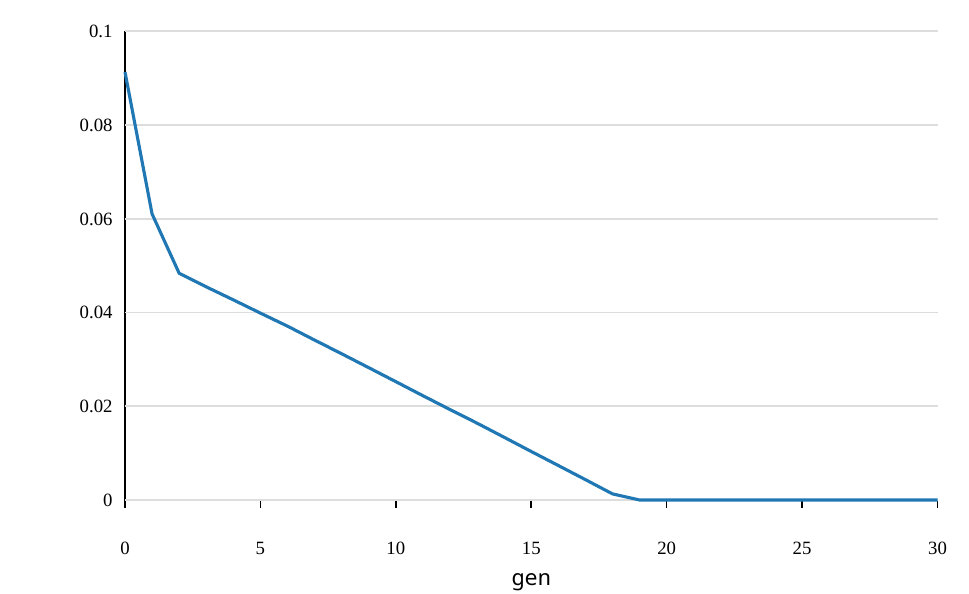}
\\[2mm]
\includegraphics[width=7cm]{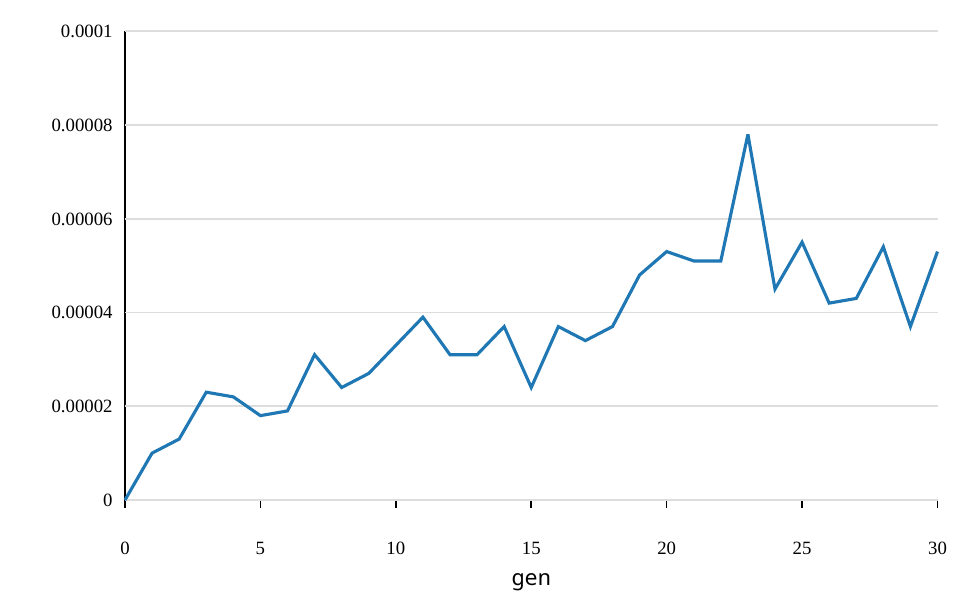}
\qquad
\includegraphics[width=7cm]{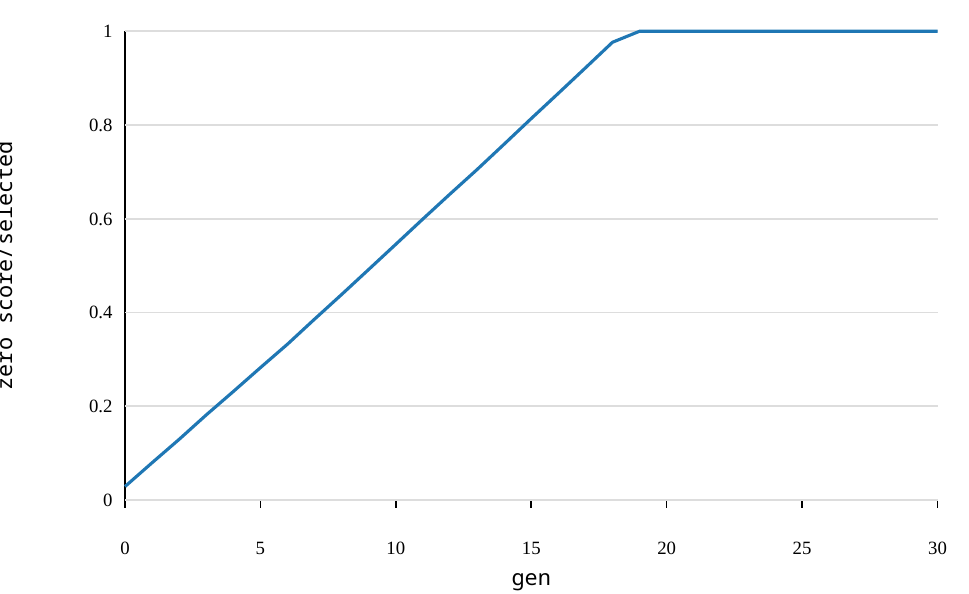}
    \caption{Loss, average score and ratio of Hadamard matrices in a typical \href{https://wandb.ai/aiformath/hadamard/runs/\runname}{$n=140$ run}.}\label{fig:lown}
\end{figure}

For moderate values of $n$ ($100$---$200$), the search for (GS type) Hadamard matrices is relatively easy, and here the contribution of PatternBoost is that once the transformer is trained, it can produce large amounts of Hadamard matrices quite fast.
See Fig.~\ref{fig:lown} for an example at $n=140$ with a \verb|sample_size| of $1\,000\,000$ and a \verb|training_size| of $50\,000$. This was run for 30 generations, but we only show 10 on the training loss plot since most of the training takes place during the first few generations.

At generation $0$ Hadamard matrices are already produced by the improvement procedure; however, the rate of production goes from $\sim 1400$ at generation 0 to $\sim 2700$ at generation 10, where it stabilises -- until the ratio of Hadamard matrices eventually reaches $1$ at generation 19 in the selected subset. Since our code picks canonical representatives in each symmetry class, we know that all produced Hadamard matrices are inequivalent; furthermore, we expect that matrices obtained in different runs are also in general distinct due to the sheer number of available ones (cf \S\ref{ssec:sym}).
E.g., out of the 12 most recent $n=140$ runs, we have 729\,732 H-matrices out of which 729\,715 are inequivalent.

The average score of the sample produced by the transformer shows an interesting effect: it goes down at first, as expected -- the minimum value, around $0.5$, is to be compared with the average value of a random sample, which is $\approx 9$, showing that the transformer has learned how to produce low score candidates -- then as the ratio of Hadamard ratio grows, it starts going up again. Presumably, the transformer starts ignoring the score statistic at some point, focusing on producing Hadamard candidates instead. Indeed, the ratio of Hadamard matrices produced by the transformer keeps going up until around generation 19 when saturation occurs. Note that this ratio remains small (of order $10^{-4}$, so about $100$ Hadamard matrices per generation).

We hope that such a ``mass generation'' of Hadamard matrices can be useful to study the typical properties of such matrices. Note, however, that since we cannot ensure uniform sampling, we do not have access to exact statistics on them.
For example, among the 729\,797 $n=140$ GS type Hadamard matrices mentioned above, the statistics on segment sums $(k_i)_{i=1,\ldots,4}$ (cf \eqref{eq:ss}; here we take the absolute value and order them so that the result is invariant under all symmetries) are:
 401\,756 $(1, 3, 7, 9)$,
 186\,827 $(1, 3, 3, 11)$,
 141\,214 $(3, 5, 5, 9)$, but it is unclear if it indicates that there are intrinsically more Hadamard matrices with segment sums $(1, 3, 7, 9)$, or if it is a bias of the transformer, or an effect of the improvement procedure. For general $n$, we observe the same tendency to generate more Hadamard matrices for which the four segment sums are spaced apart.

\subsection{Hadamard matrices for large values of $n$}
\begin{figure}[h!]
\def\runname{212_2025-11-11-15-27-50-repeat2}
    \includegraphics[width=10cm]{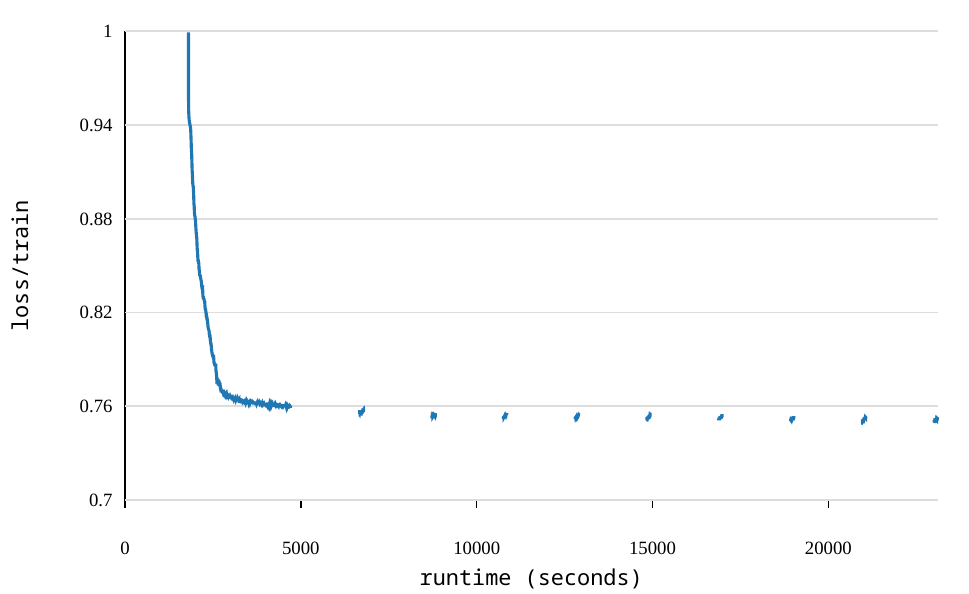}
    \\[2mm]
\includegraphics[width=7cm]{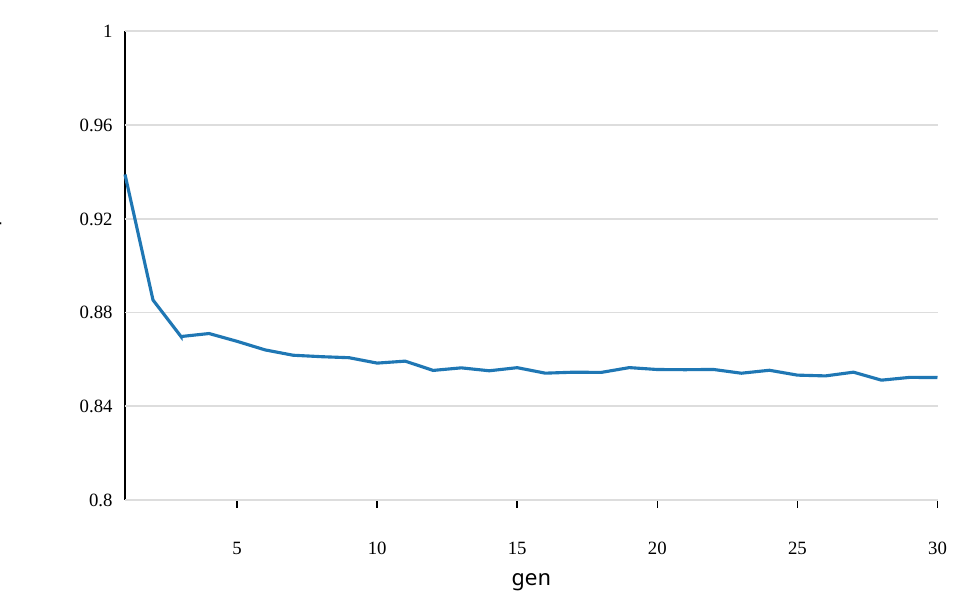}
\qquad
\includegraphics[width=7cm]{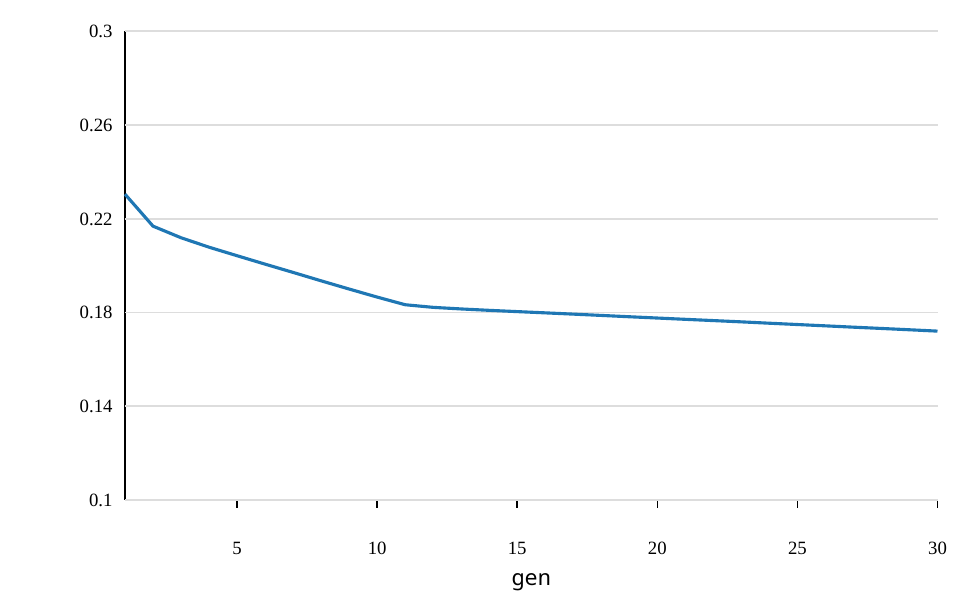}
\\[2mm]
\includegraphics[width=7cm]{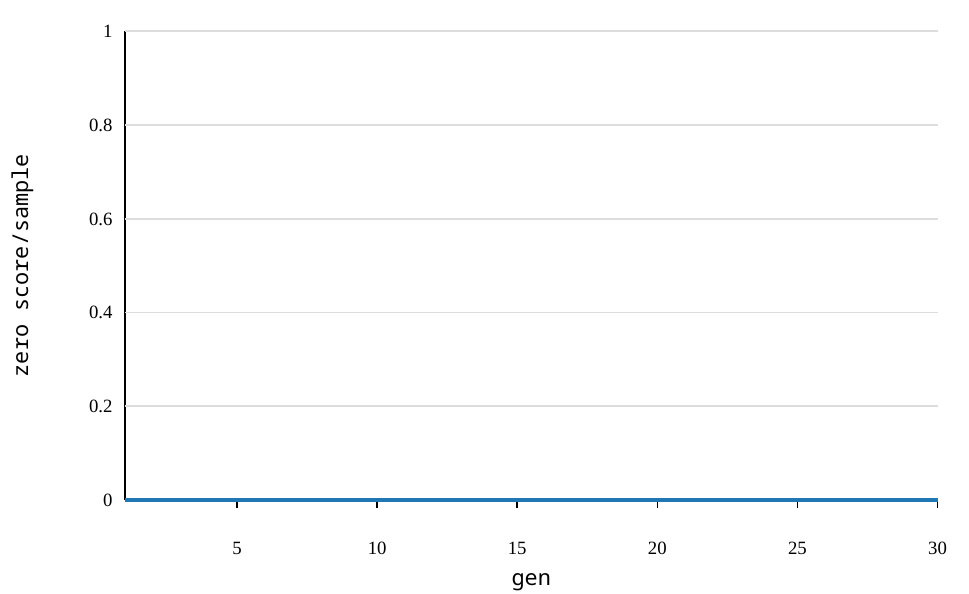}
\qquad
\includegraphics[width=7cm]{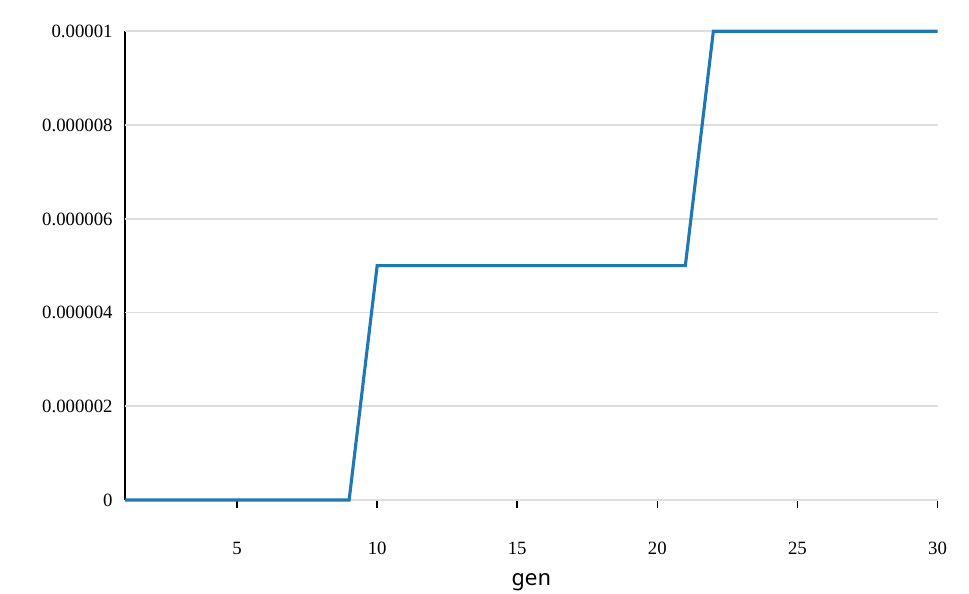}
        \caption{Loss, average score and ratio of Hadamard matrices in a typical \href{https://wandb.ai/aiformath/hadamard/runs/\runname}{$n=212$ run}.}\label{fig:hin}
\end{figure}

For values of $n$ larger than $200$ it is typically hard to produce Hadamard matrices at all. Direct exhaustive search can only be performed for very specific symmetry classes of Hadamard matrices (see, e.g., \cite{doko2024} and references therein),
but not for general GS type. We have also found that the combinatorial optimisation (of our improvement procedure) alone is not able to produce any Hadamard matrices in any reasonable amount of time, if fed random samples. However, the PatternBoost approach does succeed in producing small amounts of GS type Hadamard matrices; see Fig.~\ref{fig:hin} for a typical run at $n=212$. After six hours it produced its first Hadamard matrix; note that the transformer itself produces none, so the candidate samples it creates need to be improved to be turned into Hadamard matrices. For this range of $n$ the average score of transformer samples decreases monotonically, as expected.

Here are some samples of (the first rows of) GS type Hadamard matrices from $n=204$ to $n=252$:
\goodbreak
 
\setlength{\fboxsep}{.8pt}
\newcommand{\plusglyph}{{\fcolorbox{black}{white}{\texttt{+}}}}
\newcommand{\minusglyph}{{\fcolorbox{black}{black}{\color{white}\bfseries\texttt{-}}}}
{%
\parindent=0pt%
\parskip=0pt%
\catcode`+=13%
\catcode`-=13%
\catcode`\^^M=13%
\def^^M{\newline}%
\def+{{\plusglyph}}
\def-{{\minusglyph}}
+++-+++--++-+-+++-++--+--+----++-+----++-+++--+---+
-++-++-+++-+-+-++-+-+-+--+++---+-+------++-++++++--
+--++-+--++-+-+++++-+-+---++-----++---+--+++++++++-
--+++++-+--++--++-+-++++-++++--+-+--++++++---++++-+

--+-+---+-----+++-+-+-+--+++-----++-+-+++-++-+-+-----
----+---++---++-++--+--+----++++--++--++-----+-++-+--
--+-++------++-+-+-++--+-++-+++-++----+-++++--+--++--
--+-----+++--+++-+++---+-+--+++-+-++++++-++----+-+--+

+++++-+-+-+--+--+-++-+-+-+--+++----+---+---+----+++++--
+-----+++-+-+----+-++--+++--++--+----+++------++---+-+-
+-------++++-+-------+--+--++++-+++--++-++--++-+-++-++-
+--+-+++-+----+--+--++-+++-+--+-+-----++----+----+++--+

+++++-+-+-++-+-----++--++----+---+++---+-+-++-++++++-----
++-++-----++-+--+++-++-+++-+-+------+-+++--+++--+-++++---
--+-------++--+-+------+-+--+---+---+-+-+---+--+++++++-+-
+-+--++++-++-++-+-++-+--+--++---+-+--+++----++---++--+--+

++-+-+-++-++++++------+------++-++-++--++-++++----++-+-++-+
++-+--+++-++++---++--++--++-+-++-++--+-+---+-++--+----+-+++
++-+-+++-+++++--++-+----++++-+---+++---+-+---++--++++++---+
++++-++++-++++-+--++----+-++++-+-+-+-+--+-++-+++++--+--+++-

-+---+++--++---+-+---+-+++--++-++----+-+---++-----+-++-------
+---+-+-+-+++++--++-+--+-----++++--++-++-+++++----+-----+-+--
+-++----+-+--+---+--++-+----++-+-+--++++++--+--+++--+++--++--
+-----+++++-++-------+-+-+--+-+++-++---++-++++-+--+++-+-++-+-

++++-+-+++-+--+++-+++-++---------+++-++--++-++-++----++-----+++
---+++-+--+-++-+-+-+--+--+++-+-+++-+++---+---++-+++--+--++++-++
---+++-+-++-+--+-++++----+-----+++++--+-+-++-++-++++--++++--+++
-++++++++-+-+--+--++--+---+-+--+++-+-+-+-++++-++--+++++-+--++++
}











\subsection{Role of stacking}
One may wonder why it is not best to have the transformer produce entries of the matrix one by one, rather than the current algorithm which groups them together into ``stacks'' of $s$-bits, where $s$ is the \verb|stacking| parameter.

The answer lies in details of the transformer architecture. Recall that a decoder only transformer takes as input a sequence of tokens, and outputs probability distribution over next tokens. In the transformer architecture there is a fundamental tension between the number of tokens and the context length: with few tokens, one needs longer context length, whereas with more tokens shorter context length is possible. In the transformer the significance of tokens is encoded in their embeddings, whereas the relation between tokens is encoded in the attention mechanism. This reasoning suggests the useful heuristic that transformers will be most powerful when both the embeddings and attention are allowed to do some lifting. The heuristic is also justified by significant experimental evidence, see e.g. \cite{sennrich-etal-2016-neural} for an influential paper in the context of natural language, and \cite{nogueira2021investigating} and
\cite[\S7]{charton2022linear} for discussion in the context of AI for mathematics.

\begin{figure}[h!]
    \includegraphics[width=7cm]{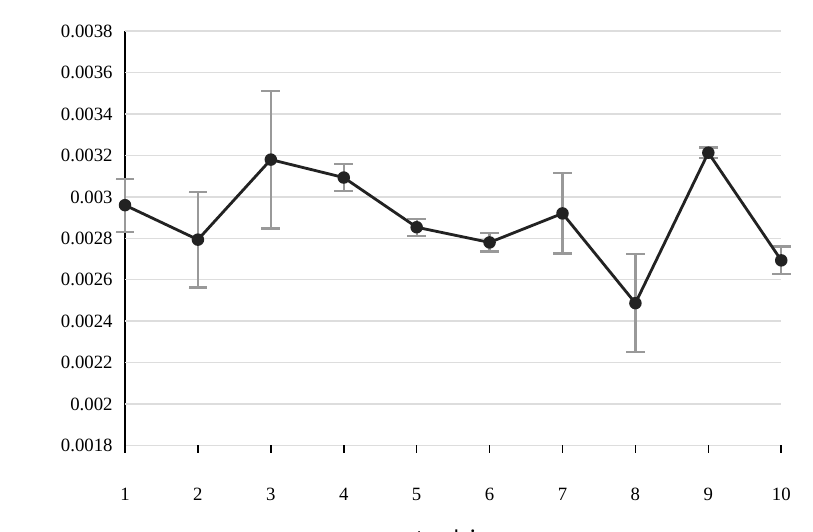}
    \qquad
    \includegraphics[width=7cm]{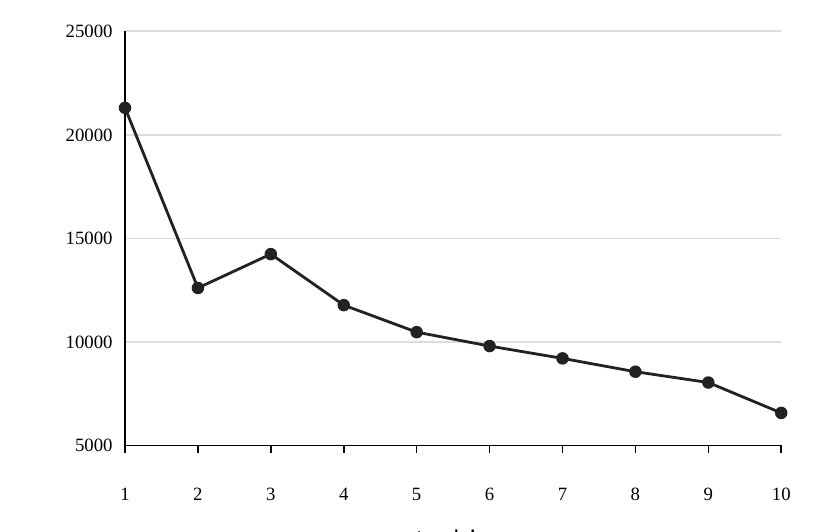}
    \caption{Influence of stacking on ratio of produced Hadamard matrices and on computation time.}
    \label{fig:stacking}
\end{figure}

Thus, we work with a larger vocabulary.
A naive estimate (minimisation of $2^s+n/s$) suggests that for our values of $n$, best results should be obtained for $s$ in the range $\{3,4,5\}$. However, there are subtle arithmetic effects due to the way the entries are packed together (i.e., related to the value of $n'\mod s$) which may affect the efficacy of the transformer. See, for example, Fig.~\ref{fig:stacking}, which shows on the left the ratio of Hadamard matrices after 30 generations for $n=172$ on a series of runs with varying $s$. There is a local maximum at $s=3$, but surprisingly $s=9$ performs equally well. Furthermore, there seems to be less variability for larger values of $s$. 

Another reason to prefer larger values of $s$ is speed; the context window of the transformer is of size $n/s$
(or $n'/s$ with the modified transformer of \S\ref{ssec:mod_transf}), which decreases with $s$, and this directly impacts the training and sampling times of the transformer.
See the right part of Fig.~\ref{fig:stacking} for the (average) run time of the same experiments at $n=172$ (using a H100 GPU).

\subsection{Learning symmetry dynamically}
Before settling on GS type matrices, we experimented with other variations of the Williamson form. For example, here is one with three circulant matrices and one anti-circulant matrix:
\begin{equation}\label{eq:alt}
M = \begin{bmatrix}
    A & B & C & DF
    \\
    -B & A & -DF & C
    \\
    -C & DF & A & -B
    \\
    -DF & -C & B & A
\end{bmatrix}
\end{equation}
The interested reader can find this on the github branch \url{https://github.com/pzinn/hadamard/tree/main-old}.
The score function $-\log\det(M/n^{1/2})$ is slightly more
complicated to write in terms of Fourier modes than \eqref{eq:score_fft}, and can be found in the code.

It is not hard to see that this Ansatz is strictly ``worse'' than the GS type, in the sense that the search space is of the same size $2^n$ but the number of Hadamard matrices of this form is strictly less than the one for GS type matrices. However one may hope that this alternate scoring function has better properties than the original one.
Though this hope turns out to be unfounded, a fascinating phenomenon appears that is worth recounting:
the transformer dynamically learns that ``most'' Hadamard matrices of this type have an additional symmetry.

\begin{figure}[h!]
    \def\runname{vhdrsqol}
    \includegraphics[width=7cm]{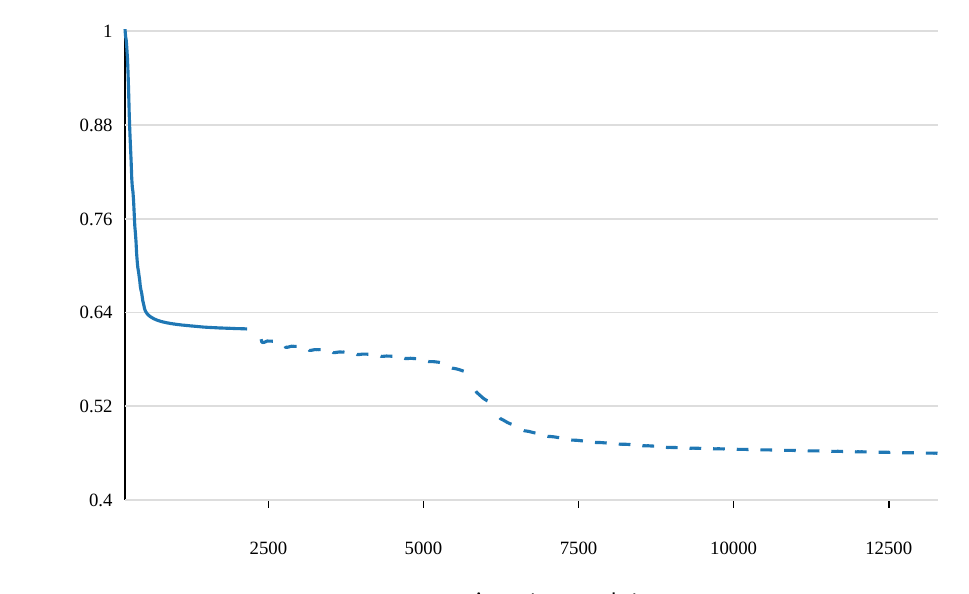}
    \qquad
    \includegraphics[width=7cm]{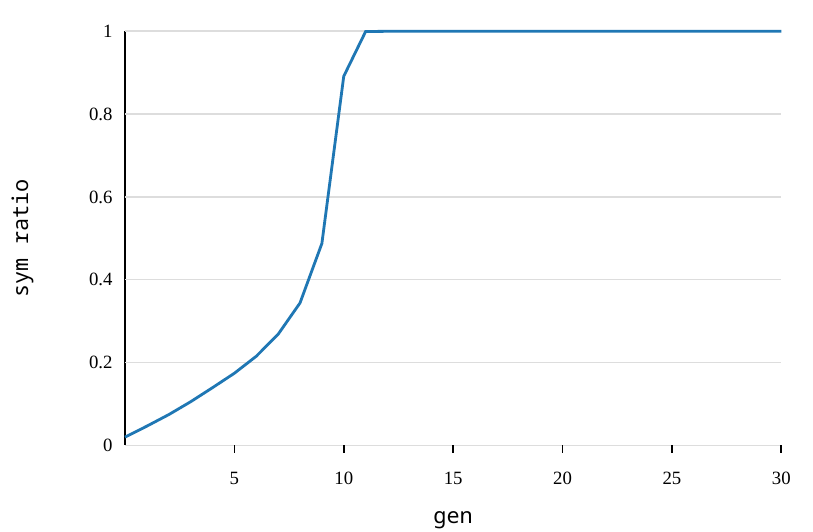}
    \\[2mm]
    \includegraphics[width=7cm]{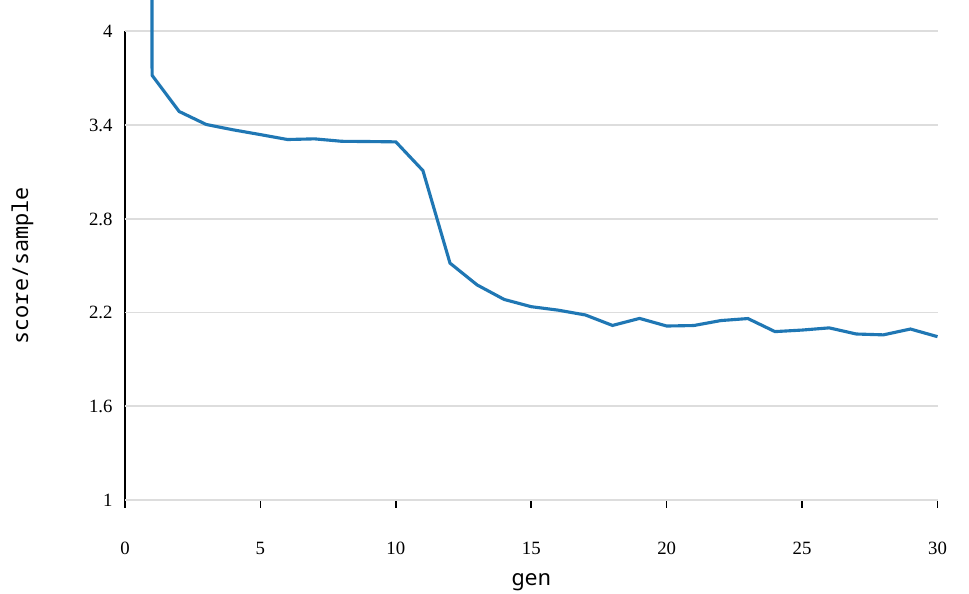}
    \qquad
    \includegraphics[width=7cm]{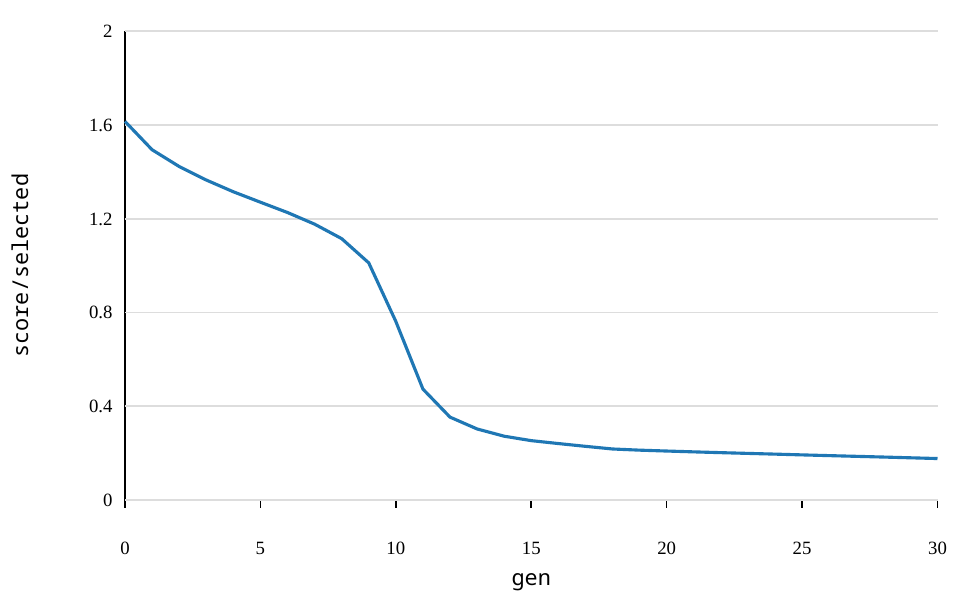}
    \\[2mm]
    \includegraphics[width=7cm]{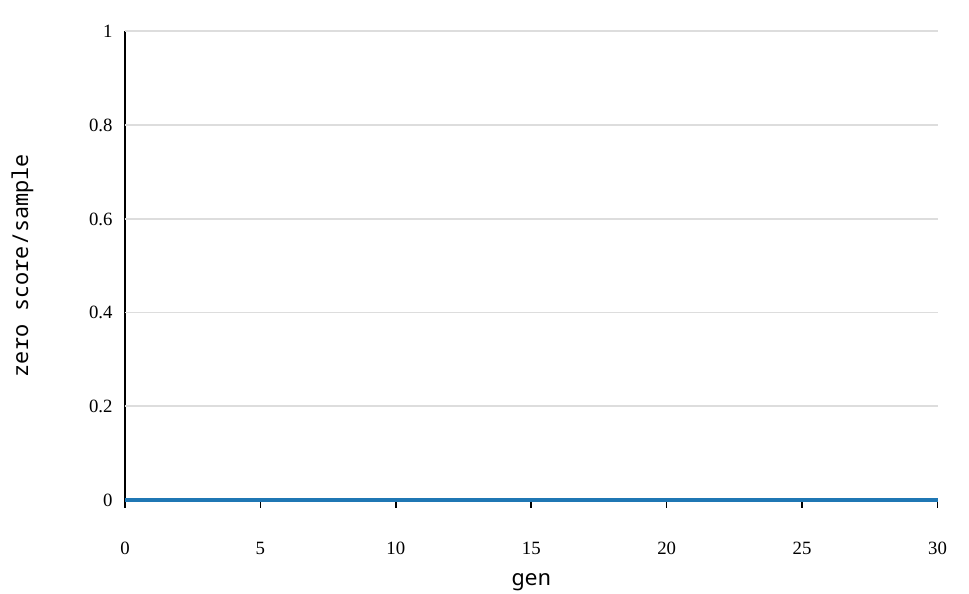}
    \qquad
    \includegraphics[width=7cm]{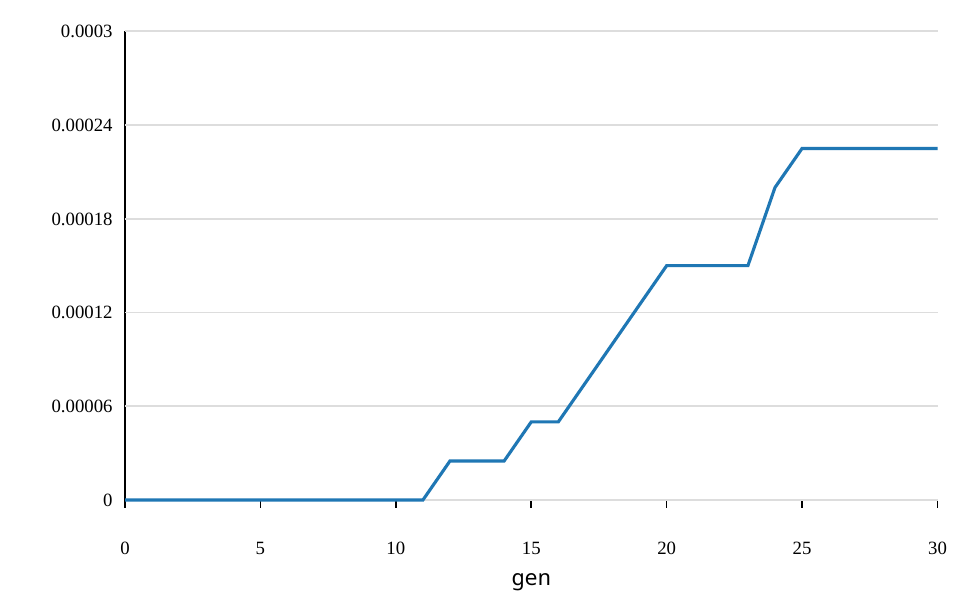}
        \caption{A typical \href{https://wandb.ai/aiformath/hadamard/runs/\runname}{$n=140$ run} with alternate Ansatz \eqref{eq:alt}.}
    \label{fig:symratio}
\end{figure}

One can see this process unfold on Fig.~\ref{fig:symratio}, which shows a typical $n=140$ run with the alternate search space of matrices of the form of \eqref{eq:alt}.
Looking at the loss, score and Hadamard ratio plots, for the first few generations, not much seems to be happening, and in particular no Hadamard matrices are produced. Then there is a sudden drop in the loss and the score, and a steady flow of Hadamard matrices appears. The explanation is given on the top right plot; it shows the ratio of selected matrices which have the following symmetry property: its first three segments are fixed by some dihedral element, i.e., they are unchanged by some reflection (the same for all three, as in this version, the automorphism group is smaller and only allows simultaneous action on the first three segments). The transformer realises at some point that it is more advantageous to produce matrices with such a symmetry in order to obtain Hadamard matrices.

At this stage, one can in principle bypass this learning phase by directly instructing the transformer to work with such symmetric matrices. This was implemented on the branch \url{https://github.com/pzinn/hadamard/tree/main-sym}. One advantage is that the search space becomes smaller -- its size is $2^{3n''+n'}$, with $n''=\lfloor n'/2\rfloor$, so roughly $2^{5n/8}$. However the results are still inferior to the standard version of the code with GS type matrices.

A related remark is that no such phenomenon occurs if we run the algorithm with GS type matrices (as on the main branch of the code), even though such symmetric matrices are also a sub-case of GS type matrices; more generally, the GS type Hadamard matrices that are produced by our code possess in general no symmetries (e.g., among the 729\,715 $n=140$ Hadamard matrices mentioned above,
only 147 have nontrivial stabiliser, of order $2$), suggesting that ``typical" GS type Hadamard matrices have no symmetries. 

\subsection{Fixed vs free segment sums}\label{ssec:discuss_ss}
As mentioned in \S\ref{ssec:ss} the code allows us to fix the ``segment sums'' of \eqref{eq:ss}. Naively, there are pros and cons to doing so. On the one hand, as already mentioned, the search space is smaller. On the other hand,
there is an increasing number of solutions to \eqref{eq:ss} as a function of $n$, so by choosing
one particular solution, one is decreasing the number of available Hadamard matrices.
Furthermore,
fixing segment sums complicates both local search (as already mentioned in \S\ref{ssec:improv}) and transformer training: concerning the latter, we have found that the transformer weights are more unstable and require smaller \verb|learning_rate|, as well as a warming up phase. This causes a longer training time.

Table~\ref{tab:ss_mod} shows the effect of varying \verb|segment_sums| for typical $n=172$ runs with a sample size of $1\,000\,000$ after 30 generations ({*} means \verb|segment_sums| is turned off, which leaves the segment sums unconstrained). The result is unambiguous: the number of produced Hadamard matrices is significantly higher with the fixed segment sums.
Note that $(1,5,5,11)$ (resp.\ $(5,7,7,7)$) is the most (resp.\ least) frequent occurrence among segment sums in experiments with \verb|segment_sums| turned off.


\begin{table}[h]
\begin{tabular}{c|r|r|r}
& {\quad*\quad} & {(5,7,7,7)} & {(1,5,5,11)} \\\hline
\multirow{2}{*}{\ttfamily False} & 1 & 0 & 1 \\
                                & 270 & 1\,009 & 1\,087 \\\hline
\multirow{2}{*}{\ttfamily True}  & 14 & 17 & 26 \\
                                & 305 & 1\,196 & 1\,437
\end{tabular}
\bigskip
    \caption{Influence of {\tt transformer\_uses\_score} (rows) and {\tt segment\_sums} (columns) on the generation of $n=172$ Hadamard matrices:
    top/bottom number is pre/post improvement total number of Hadamard matrices after 30 generations.}
    \label{tab:ss_mod}
\end{table}

\subsection{Comparison of vanilla and modified transformers}\label{ssec:discuss_mod}
The same Table~\ref{tab:ss_mod} also shows the effect of using either the vanilla transformer or the modified one, see \S\ref{ssec:mod_transf} for details.
It is striking that at $n=172$, the vanilla transformer barely produces any Hadamard
matrices at all, whereas the modified one produces a few (of the order of 1 per generation). After improvement the numbers of Hadamard matrices are typically increased by about 20\%.

A slightly 
counterintuitive feature is that the average score of transformer samples is significantly worse with \verb|transformer_use_score| turned on.
One possible interpretation is that the mismatch between the training (which typically uses data with positive score) and the sampling (which tries to generate zero score samples) means that the modified transformer tries to extrapolate rather than simply sample from the known population, so is more prone to riskier choices (but potentially more rewarding ones too).
This deserves further investigation.

\bibliographystyle{alpha}
\bibliography{gen}

\end{document}